\documentclass[12pt]{article}
\usepackage{amssymb, amscd}
\newcommand{\Pe}{\mathbb{P}}

\newcommand{\Q}{\mathbb{Q}}

\newcommand{\co}{{\mathcal{O}}}

\newcommand{\I}{{\mathcal{I}}}
\newcommand{\ra}{\rightarrow}

\newcommand{\vp}{\varphi}

\newcommand{\id}{\mbox{\rm id}}

\newcommand{\Pic}{\mbox{Pic}}

\begin{document}
\def\flexleftarrow#1{\mathtop{\makebox[#1]{\leftarrowfill}}}
\begin{center}
{\Large \bf Severi varieties  and branch curves }
\end{center}
\begin{center}
 {\Large \bf of abelian surfaces  of type $(1,3)$}\\
\end{center}
\bigskip
\begin{center}
H. Lange  and E. Sernesi \footnote{Research supported by the joint Vigoni Program
of DAAD and CRUI}
\end{center}
\vskip2cm
\begin{center}
{\bf Introduction}
\end{center}
Let $(A,L)$ be a principally polarized abelian surface of type
$(d_1,d_2)$. The linear system $|L|$ defines a map
$\vp_L:A \to \Pe^\ast_{d_1d_2-1}= \Pe(H^0(L)^\ast)$.
The geometry of this map is well understood except for 
$(d_1,d_2)=(1,3)$ (see [CAV]). In this case $\vp_L$ is a $6:1$ covering of
$\Pe_2^\ast$ with branch curve $B\subset \Pe_2^\ast$ of degree 18. 
The first problem is to study the curve $B$.  
The main result of this paper  is that for a general abelian surface 
$(A,L)$ of type $(1,3)$ the curve $B$ is irreducible, admits 72 cusps,
36 nodes or tacnodes, each tacnode counting as two nodes, 72 flexes and
36 bitangents. The main idea of the proof is to use the fact 
that for a general $(A,L)$ of type $(1,3)$
the discriminant curve $V \subset \Pe_2 = |L| = \Pe(H^0(L))$, coincides
with the closure
of the Severi variety $V_{L,1}$ of curves in $|L|$ admitting a node, 
and is dual to the curve $B$ in the sense of projective geometry. We 
investigate $V$ and $B$ via degeneration to a special abelian surface. 
\medskip

In more detail the contents of the paper is as follows.
In section 1  some well known facts about families of nodal curves on
surfaces are recalled and proved for the case of an abelian surface. 
In section 2 we consider curves of genus 2 in a linear system $|L|$ as
above. There are only finitely many such curves and we prove that they
 have only ordinary singularities. 
In section 3 we introduce and study the incidence curve.  We explain
how it is  related to $B$ and $V$ and  compute its numerical characters. 
Section 4 is devoted to a detailed analysis of a special case: we consider
the abelian surface $A=E\times E$, where $E$ is an elliptic curve, suitably
polarized by a line bundle of type $(1,3)$, and we describe completely
the curves $B$ and $V$. In the final section 5 we take up the general case:
we degenerate a general $(A,L)$ of type $(1,3)$ to the surface 
considered in \S 4 and, using the analysis made in that special case, we 
prove our main result. 
\medskip

We are grateful to  Mike Roth for suggesting the use of the 
incidence curve  $\Gamma$ in \S 3, and for some helpful conversation. 

\newpage
\begin{center}
{\bf \S 1. Generalities on Severi varieties}
\end{center}
We recall a few  known facts on families of nodal curves 
on an algebraic surface.
 
Let $S$ be a nonsingular projective connected algebraic surface,
$L\in \Pic(S)$ such that the complete linear system $|L|$ has a
nonsingular and connected general member. The linear system 
$|L|$ is contained in a larger system of divisors, the 
{\it continuous system} defined by $L$, which 
 will be denoted by  $\{L\}$ in  what follows. By definition 
 $\{L\}$ is the connected component of Hilb$^S$ containing $|L|$. 
 If
$X\subset S$ is a curve belonging to the system  $\{L\}$ we will 
denote by $[X]$ the point of $\{L\}$ parametrizing $X$.

 For each integer $\delta \geq 0$ denote by $V_{L,\delta} \subset |L|$
 the functorially defined locally closed subscheme which parametrizes 
 the family of all curves in $|L|$ having $\delta$ nodes 
 and no other singularity (see [W]); they will be 
 called $\delta$-{\it nodal curves} . We likewise denote by 
 $V_{\{L\},\delta} \subset \{L\}$ the analogous closed subscheme  
 of $\{L\}$ parametrizing $\delta$-nodal curves. 
 
 It is customary to call {\it regularity} the condition of being
 nonsingular of codimension $\delta$ for 
 $V_{L,\delta}$ or $V_{\{L\},\delta}$. 

Let $X \in V_{L,\delta}$ and let $N\subset S$ be the  scheme of its nodes,
which of course consists of $\delta$ distinct reduced points.
We have an exact commutative diagram of sheaves on $S$:
$$
\matrix{&0&&0&&& \cr
&\downarrow&&\downarrow&&& \cr
&\co_S&=&\co_S&&& \cr
&{\sigma_X}\downarrow\\ &&\\ \downarrow{\sigma_X}&&& \cr
0 \to &{\cal I}_{N/S}\otimes L& \to & L & \to & L\otimes\co_N & \to 0 \cr
& \downarrow&&\downarrow&& \Vert\wr&& \cr
0 \to &{\cal I}_{N/X}\otimes N_X& \to & N_X & \to &   T^1_X & \to 0 \cr
&\downarrow&&\downarrow&&& \cr
&0&&0&&&}
$$
where $N_X = L\otimes\co_X$ is the normal sheaf of $X$ in $S$, and
$ T^1_X$ is the first cotangent sheaf of $X$ (see [W]).
Consider the following diagram obtained by taking cohomology of it:
$$
\matrix{
0\to & {H^0(\I_{N/S}\otimes L) \over \langle\sigma_X\rangle} & \to &
{H^0(L) \over \langle\sigma_X\rangle} & 
\buildrel h \over \longrightarrow & H^0(L\otimes\co_N) \cr\cr
& \bigcap&& \bigcap && \Vert \cr\cr
0\to & H^0(\I_{N/X}\otimes N_X) & \to & H^0(N_X) & \to & H^0(T^1_X) \cr\cr
&\downarrow&&\downarrow&& \cr\cr
&H^1(\co_S)&=&H^1(\co_S)&& \cr\cr
&\downarrow&&\downarrow&& \cr\cr
  & H^1(\I_{N/S}\otimes L) &   &0&&\cr\cr
&\downarrow&&&&\cr\cr
&H^1(\I_{N/X}\otimes N_X)&&&&\cr\cr
&\downarrow&&&&\cr\cr
&H^2(\co_S)&&&&\cr\cr
&\downarrow&&&&\cr
&0&&&&}
$$
where the zero at the left bottom is because $H^2(\I_{N/S}\otimes L)=0$.
We have identifications:
$$
{H^0(\I_{N/S}\otimes L) \over \langle\sigma_X\rangle} = 
{\bf T}_{[X]}(V_{L,\delta})
$$  
$$
{H^0(L) \over \langle\sigma_X\rangle} = {\bf T}_{[X]}(|L|)
$$
$$
H^0(\I_{N/X}\otimes N_X) = {\bf T}_{[X]}(V_{\{L\},\delta})
$$ 
and
$$
H^0(N_X) = {\bf T}_{[X]}(\{L\})
$$
where  ${\bf T}_{[X]}(-)$  denotes the tangent space of \ $-$ \ at $[X]$.
Therefore the above diagram is read as:
$$
\matrix{
0\to & {\bf T}_{[X]}(V_{L,\delta}) & \to &
{\bf T}_{[X]}(|L|) & 
\buildrel h \over \longrightarrow & H^0(L\otimes\co_N) \cr\cr
& \bigcap&& \bigcap && \Vert \cr\cr
0\to & {\bf T}_{[X]}(V_{\{L\},\delta}) & \to & {\bf T}_{[X]}(\{L\}) & \to 
& H^0( T^1_X)\cr\cr
&\downarrow&&\downarrow&& \cr\cr
&H^1(\co_S)&=&H^1(\co_S)&& \cr\cr
&\downarrow&&\downarrow&& \cr\cr
  & H^1(\I_{N/S}\otimes L) &   &0&&\cr\cr
&\downarrow&&&&\cr\cr
&H^1(\I_{N/X}\otimes N_X)&&&&\cr\cr
&\downarrow&&&&\cr\cr
&H^2(\co_S)&&&&\cr\cr
&\downarrow&&&&\cr
&0&&&& } \eqno (1)
$$
and all the maps between tangent spaces are the differentials 
of the corresponding inclusions.
Since we have codim$_{|L|}(V_{L,\delta}) \le \delta$, 
we see that 
 
\bigskip\noindent
{\bf Proposition 1.1.} 
$$
\hbox{$V_{L,\delta}$ is regular at $[X]$} 
\Leftrightarrow \hbox{$h$ is surjective}
$$
$$
\Leftrightarrow 
\hbox{ $N$ imposes independent conditions to $|L|$}
\Leftrightarrow H^1(\I_{N/S}\otimes L)=0
$$
\medskip 

Note that if $S$ is a regular surface then 
$V_{L,\delta}=V_{\{L\},\delta}$ and the vertical inclusions are 
equalities.

\bigskip

Let $A$ be an abelian surface and $L$ an ample line bundle of type $(d_1, d_2)$ on $A$. We call $d :=
d_1 d_2$ the {\it degree} of $L$. According to Riemann-Roch the linear system $|L|$ is of dimension
$d-1$, which we will assume to be $\geq 2$. 
We always assume that $|L|$ has no fixed component. This means that the polarized abelian variety
$(A,L)$ is not of the form $(E_1 \times E_2, p^{\ast}_1 L_1 \otimes p^{\ast}_2 L_2)$ with elliptic
curves $E_1$ and $E_2$ and  line bundles $L_1$ of degree 1 on $E_1$ and $L_2$  of degree $d$ on $E_2$.
According to [CAV] Proposition 10.1.3 the general member of $|L|$ is a smooth irreducible curve on $A$ of
genus $d+1$. 
We have
$$
{\rm dim}(\{L\}) = {\rm dim}(|L|)+2 = d+1
$$
 because  $\{L\}$ 
consists of a two-dimensional family of 
translates of $|L|$. More precisely, denoting by $\hat A$ the dual of $A$, 
we have a surjective morphism:
$$
\psi: \{L\} \to \hat A
$$ 
whose fibres are linear systems. 
For the same reason
$$
 {\rm dim}(V_{\{L\},\delta}) = {\rm dim}(V_{L,\delta}) + 2
$$ 
 The curve $X$ has arithmetic genus 
 $$
 p_a(X) = {1\over 2}L^2 +1= d+1
 $$
 and $N_X = \omega_X$ is its dualizing sheaf. Therefore 
 $$
 h^0(N_X) = p_a(X) = d+1 = dim(\{L\})
 $$ 
 and it follows that $\{L\}$ is nonsingular of dimension $d+1$ at  $X$.

 The differential of $\psi$ at $X$ is the coboundary map
 $$
 H^0(N_X) \to H^1(\co_A)
 $$
 which is surjective because $H^1(L)=0$. It follows that $\psi$ is
 smooth at $[X]$ with fibre $\psi^{-1}(\psi(X)) = |L|$.

 In order to compute  $ {\bf T}_{[X]}(V_{\{L\},\delta})$ consider the normalization
map $\nu: C \to X$. Then we have an isomorphism
$$
H^0(X,\I_{N/X}\otimes N_X) = H^0(X,\I_{N/X}\otimes \omega_X)
\cong H^0(C,\omega_C)
$$
Then
$$
{\rm dim} \big[{\bf T}_{[X]}(V_{\{L\},\delta})\big] = h^0(C,\omega_C) = g
$$
where $g$ is  the geometric genus of $C$, which is given by the formula
$$
g = p_a(X) - \delta + c -1
$$
where $c$ is the number of connected components of $C$. Since 
$$
{\rm dim}\big[V_{\{L\},\delta}\big] \ge p_a(X) - \delta
$$
 we have:

\bigskip\noindent
{\bf Proposition 1.2.} {\it  Let $(A,L)$ be as above. Then
$V_{\{L\},\delta}$ is regular
at $[X]$ if and only if $X$ is irreducible.}
\par
 
\medskip
In particular we see that the conditions for regularity of 
$V_{L,\delta}$ (Prop. 1.1) and of $V_{\{L\},\delta}$
 (Prop. 1.2) are different. Nevertheless we have the following:

\bigskip\noindent
{\bf Proposition 1.3}  {\it In the same situation of Prop. 1.2, assume that
$X$ is reducible. Then $V_{L,\delta}$ is not regular at $[X]$.}
 
 \medskip
\noindent
{\it Proof.}  Diagram (1)  shows that if $h^1(\I_{N/X}\otimes N_X) \geq 2$ then 
$h^1(\I_{N/S}\otimes L) \neq 0$ and $V_{L,\delta}$ is not regular by
Prop. 1.1. Since from the above analysis  it follows that 
 $h^1(\I_{N/X}\otimes N_X)$ is equal to the number of irreducible components of $X$, the conclusion follows.
\hfill $\square$

 \medskip
Summarizing, irreducibility of $[X]$ is a necessary and sufficient condition for the regularity of $V_{\{L\},\delta}$, while it is only a necessary condition for the regularity of  $V_{L,\delta}$.

\medskip\noindent
 {\bf Remark.} 
 Note that  Propositions 1.1 and 1.2 apply as well to the case when
 $S$ is a K3-surface. In this case 
 $V_{L,\delta}=V_{\{L\},\delta}$ and we obtain the following 
 combination of the two propositions:

 \medskip\noindent
 {\bf Proposition 1.4.} {\it Let $S$ be a K3-surface, $L$ a 
 globally generated line bundle on $S$ such that dim$(|L|)\geq 2$,
 and $X\in V_{L,\delta}$. 
 Then $V_{L,\delta}$ is regular at $[X]$ if and only if 
 $X$ is irreducible. }
 
\medskip
 In the sequel we   will be especially interested in  $V_{L,1}$
and in its closure $\overline{V}_{L,1}$.  
From the assumptions made on
$L$ it follows that  $|L|$ is base point free and therefore 
Proposition 1.1 implies that  $V_{L,1}$ is regular (i.e.
nonsingular of codimension 1) everywhere. 

We will also consider the {\it discriminant locus} $V \subset |L|$
(see \S 3 for its scheme-theoretic definition)
which parametrizes all singular curves of $|L|$.  Of course we have
$\overline{V}_{L,1}\subset V$ but the inclusion is proper  in general
(see \S 4). 

%\newpage
\vspace{2cm}
\begin{center}
{\bf \S 2 Singularities of curves of genus 2 on abelian surfaces.}
\end{center}
Let $(A,L)$ be an abelian surface of type  $(d_1,d_2)$. In this section we show that every  curve  of genus 2 in $|L|$ admits at most ordinary singularities. 
For this we need
the following result on principally polarized abelian surfaces which is certainly well-known. Not
knowing a reference we give a proof.
\medskip

\noindent
{\bf Proposition 2.1.} {\it Let $(J, \co_J(C))$ denote a principally polarized abelian surface. If
translates $t^{\ast}_x C$ and $t^{\ast}_y C$ of the curve $C$ do not have a component in common, they
intersect exactly in 2 distinct points.}
\medskip

\noindent
{\it Proof.} $J$ is either a product of 2 elliptic curves $E_1 \times E_2$ and $C = E_1 \times \{ 0 \} +
\{ 0 \} \times E_2$ or the Jacobian of a smooth curve $C$ of genus $2$. The assertion being obvious in
the product case, we may assume that $C$ is smooth of genus 2 containing $0 \in J$. Translating and
passing to an algebraically equivalent line bundle we see that it suffices to show that $\# (C \cap
t^{\ast}_y C) = 2$  only in the case $0 \in C \cap t^{\ast}_y C$ and $C \not= t^{\ast}_y C$.
\smallskip

\noindent
Note first that the dual abelian variety $\widehat{J}$ parameterizes the translates of $C$: 
$\widehat{J} = \{ [ t^{\ast}_y C] \, | \, y \in J \}$. In particular $t^{\ast}_y C \not= C$ for all $y
\not= 0$. Moreover, $C^2 = 2$ means that $C \cap t^{\ast}_y C$ consists of 2 points
counted with multiplicities.
\smallskip

\noindent
Now consider the closure $\widetilde{C}$ of the set $\{ (C \cap t^{\ast}_y C) \backslash \{ 0 \} \, | \,
0 \in C \cap t^{\ast}_y C \}$. $\widetilde{C}$ is certainly of dimension 1 contained in $C$ and as $C$
is irreducible, we have $\widetilde{C} = C$. Note that
$$
0 \in C \cap t^{\ast}_y C \Longleftrightarrow y \in C .
$$
Hence the map $C - \{ 0 \} \ra \widetilde{C} , \,\, y \mapsto (C \cap t^{\ast}_y C) \backslash \{ 0 \}$
extends to a morphism
$$
\vp : C \ra \widetilde{C} = C.
$$
But $\vp$ is not constant. So by Hurwitz' formula it has to be an isomorphism. This means that $(C \cap
t^{\ast}_y C \backslash \{ 0 \}) \not= 0$ if $y \not= 0$, which was to be shown. \hfill $\square$
\medskip

\noindent
For the rest of this section let $L$ denote an ample line bundle of type $(d_1, d_2)$ on the abelian surface $A$. By a 
{\it curve of
genus} 2 we mean an irreducible reduced curve of geometric genus 2. As a consequence of Proposition 2.1
we obtain:
\medskip

\noindent
{\bf Proposition 2.2.} {\it Any curve of genus $2$ in the linear system $|L|$ admits at most ordinary
singularities of multiplicity $\leq \frac{1}{2} (1 + \sqrt{8 d_1 d_2 - 7})$.}
\medskip

\noindent
{\it Proof.} Suppose $\overline{C} \subset |L|$ is of genus 2 and $C \ra \overline{C}$ its
normalization. The universal property of the Jacobian yields an isogeny $f: J = J(C) \ra A$ such that
after a suitable embedding $C \hookrightarrow J$ the following diagram commutes
$$
\begin{array}{ccc}
C & \subset & J\\
\downarrow && \downarrow {\scriptstyle f}\\
\overline{C} & \subset& A
\end{array}
$$
$f$ being \'{e}tale, it is clear that $\overline{C}$ cannot admit a cusp. Suppose $\overline{C}$ admit a
tacnode, i.e. a point $x$ such that $\overline{C}$ has 2  linear branches $\overline{C}_1$ and $\overline{C}_2$
touching  in $x$.
\smallskip

\noindent
Let $x_1$ and $x_2$ be 2  preimages of $x$ in $C$. The branches $C_1$ and $C_2$ of $C$ in $x_1$ and
$x_2$ map to $\overline{C}_1$, and $\overline{C}_2$. Since $t^{\ast}_{x_1 - x_2} (x_1) = x_2$ this
implies that $t^{\ast}_{x_1 - x_2} C$ and $C$ intersect in $x_2$ of multiplicity $2$. But this
contradicts Proposition 2.1. Hence $\overline{C}$ admits only ordinary 
singularities. Suppose $\overline{C}$ admits an ordinary singularity of multiplicity $\nu$.
Then
$$
2 \leq g (\overline{C}) = d_1 d_2 +1 - \frac{\nu (\nu -1)}{2}
$$
since there is no curve of geometric genus $\leq 1$ in $|L|$. But 
$\frac{1}{2} (1 \pm \sqrt{8 d_1 d_2 -7)}$ are
the roots of the polynomial $\nu^2 - \nu  - 2 d_1 d_2 + 2 $ in $\nu$, implying $\nu \leq \frac{1}{2} ( 1 +
\sqrt{8 d_1 d_2 -7})$. \hfill $\square$
\medskip

\noindent
{\bf Remark 2.3.}  In the situation of Prop. 2.2 there are only finitely many
curves $\overline{C}$ of geometric genus 2 in $|L|$. This depends on the fact that each such curve defines an isogeny of bounded degree
from the Jacobian of the normalization $C$ of $\overline{C}$ to $A$ and
it is defined by such an isogeny up to translation.
Since  there are only finitely many such isogenies and also finitely many 
translations of $\overline{C}$ which are still in $|L|$, the conclusion
follows.

%\newpage
\vspace{2cm}
\begin{center}
{\bf \S 3 The incidence curve.      }
\end{center}
Let $(A,L)$ be an abelian surface of type  $(d_1,d_2)$,
where as usual we assume
that  $L$ defines an irreducible polarization, i. e. the linear system $|L|$ has no fixed
components. 

\noindent
Let $\Pe_2\subset |L|$ be a general net in $|L|$.
The {\it incidence curve} $\Gamma$ associated to  $\Pe_2$ is defined  as follows.

Let $J_1 (L)$ denote the first jet bundle of $L$. The fibre of $J_1(L)$ at a point $x \in A$ is the
space
$L \otimes \co_A / I^2_x$ of 1-jets of sections of $L$. It fits into the basic exact sequence
$$
0 \ra \Omega^1_A \otimes L \ra J_1 (L) \ra L \ra 0. \eqno(1) 
$$
There is a natural homomorphism of sheaves
$$
\sigma : \pi^{\ast}_2 \co_{\Pe_2} (-1) \ra \pi^{\ast}_1 J_1 (L)
$$
defined by associating to every local section the constant and linear terms of its Taylor expansion.
This gives a global section of the locally free sheaf
$$
E = \pi^{\ast}_2 \co_{\Pe_2} (1) \otimes \pi^{\ast}_1 J_1 (L).
$$
 By definition $\Gamma$ is the  vanishing scheme of this section and 
the third Chern class of $E$ is the  class of
 $\Gamma$. Note that every component of $\Gamma$ has dimension at most 1
since there are at most finitely many nonreduced curves in
$\Pe_2$. On the other hand, being the zero locus of a section of a rank 3
vector bundle on $A\times \Pe_2$, it has pure dimension 1 and is a local 
complete intersection curve. 
Set theoretically
$$
\Gamma := \{ (x, [C]) \in A \times \Pe_2 \, | \, x   \quad \mbox{is a singular point of} \quad C \}.
$$

\medskip\noindent
Let $\pi_1$ and $\pi_2$ denote the natural projections of $A \times \Pe_2$.  
Consider  the  map $\vp_{\Pe_2} : A \ra \Pe^{\ast}_2$ associated to $\Pe_2$,
which  is a covering of degree $2d$. Here $ \Pe^{\ast}_2$ denotes the 
dual projective space.

\medskip\noindent
{\bf Lemma 3.1.} {\it
 $\pi_1(\Gamma)$ is equal to the ramification locus $R$
of $\vp_{\Pe_2}$.}

\medskip\noindent
{\bf Proof.}  If $y \in R$ then $\Pe_2(-y)$ (the pencil in $\Pe_2$
of curves passing through $y$) consists of curves with a fixed tangent at $y$:
therefore  $\Pe_2(-y)$ contains a curve $C$ singular at $y$ because such a curve 
satisfies only one extra linear condition. Therefore $(y,[C])\in \Gamma$
and $\pi_1(y,[C])=y$. \hfill $\square$

\medskip\noindent
Denoting $l := c_1 (\pi^{\ast}_1 (L))$ and $h := c_1 (\pi^{\ast}_2 \co_{\Pe_2} (1))$ we have
 $$
 l^2 h^2 = 2d,
 $$
 whereas all other intersection numbers vanish. Using (1) we get
 $$
 c_1 (E) = 3(l+h) \quad \mbox{and} \quad \Gamma  = c_3 (E) = 3  (l^2 h + l h^2).
 $$
 There is an isomorphism
 $$
   \omega_{\Gamma} = \omega_{A \times \Pe_2}\otimes \bigwedge^3 N_{\Gamma}  \eqno(2)
 $$
 Now 
 $$
 c_1 (\omega_{A \times \Pe_2} | \Gamma) = -3 h \cdot \Gamma = -9 l^2 h^2 = -18 d.
 $$
 On the other hand $N_{\Gamma} = E|\Gamma$, since $\Gamma$ is the top Chern class of the vector bundle
 $E$. This gives
 $$
 c_1 (N_{\Gamma}) = c_1 (E) \cdot \Gamma = 9 (l+h) (l^2 h + l h^2) = 36 d .
 $$
 So (2) yields
 $$
 \deg \omega_{\Gamma}  =  18 d
 $$
and we have proven
\medskip

\noindent
{\bf Proposition 3.2.} $p_a (\Gamma) = 9 d+1.$
\medskip

\noindent
The equality 
$$
\pi_2 (\Gamma)=V
$$
defines scheme theoretically  the {\it discriminant locus $V$ of the net}. 
Since    $\Gamma$ 
is a connected curve, we have   the following:

\medskip\noindent
{\bf Lemma 3.3.} {\it $V$ is purely one-dimensional, i.e. it is a plane curve.}

\medskip

 On the other hand consider the map $\vp = \vp_{\Pe_2} : A \ra \Pe^{\ast}_2$
associated to the net $\Pe_2$. 
As noted
above, $\vp$ is a covering of degree $2d$. Restricting to a line in $\Pe^{\ast}_2$ and using Hurwitz'
formula, one sees that the branch divisor $B \subset \Pe^{\ast}_2$ of $\vp$ is a curve of degree $ 6d$. 

\medskip
From now on we will restrict to the case $(d_1,d_2)=(1,3)$. Therefore $|L|$ is 
a net and the branch curve $B$ has degree 18.
Let's consider the curve $V$. We have the following

\medskip\noindent
{\bf Proposition 3.4.}  {\it 
$\overline{V}_{L,1}$ has degree $\leq 18$ 
and equality holds if and only if   $V = \overline{V}_{L,1}$. In this case 
$\pi_2:\Gamma \ra V$ is birational on each irreducible component of $\Gamma$.
In particular, if $\Gamma$ is reduced and $V = \overline{V}_{L,1}$
then $V$ is also reduced.}

\medskip\noindent
{\it Proof.} 
Let $\Pe_1 \subset |L|$ be a general pencil.  The pencil $\Pe_1$ being general, it
has $6$ base points $x_1, \ldots , x_6$. No curve in the pencil is singular at one of the base points,
since otherwise $L^2 = C_1 \cdot C_2 > 6 $. Let $M$ denote the blow-up of $A$ in the $6$ base points.
$f : M \ra \Pe_1$ is a fibration with smooth general fibre $C_{\rm gen}$. For any $s \in \Pe_1$ we denote
$C_s := f^{-1} (s)$, the fibre over $s$.    We have the following
relation between the topological Euler characteristics
$$
e(M) = e(C_{\rm gen}) \cdot e(\Pe_1) + \sum_{s \in \Pe_1} (e(C_s) - e(C_{\rm gen})) .
$$
(see [BPV] Proposition III.11.4). But $\sum_{s \in \Pe_1} (e(C_s) - e(C_{\rm gen}))$ is the number of singular curves in the pencil, each 
counted with some multiplicity, and the nodal curves 
have multiplicity 1. Hence we get
$$
 \sum_{s \in \Pe_1} (e(C_s) - e(C_{\rm gen})) 
\geq \deg \overline{V}_{L,1} \eqno (3)
$$
Moreover 
$$
\sum_{s \in \Pe_1} (e(C_s) - e(C_{\rm gen})) 
= 6 - (2 - 2(4)) \cdot 2 = 18 
$$
since $C_{\rm gen}$ is smooth of genus $4$. 
The  inequality in (3) is an equality if and only if all curves 
in the pencil are 1-nodal and this means that 
$V = \overline{V}_{L,1}$. In this case a general element of any irreducible
component of $V$ is 1-nodal, and therefore $\pi_2$ is birational on each
component of $\Gamma$. \hfill $\square$

%\medskip\noindent
%{\bfRemark 3.5.} Actually the above proof shows more. Infact the only 
%singular curves in 
%the pencil which contribute with multiplicity 1 are the 1-nodal ones. 
%Therefore if a component of 

%\newpage
\vspace{2cm}
\begin{center}
{\bf \S 4 Abelian surfaces of type \boldmath$(1,3)$\unboldmath, a special case.}
\end{center}
%\end{document}
Let $(A,L)$ be an abelian surface of type $(1,3)$, so the linear system $|L|$ is a plane $\Pe_2 = P(H^0 (L))$. According to \S 3 
$\vp = \vp_L : A \ra \Pe^{\ast}_2 = P(H^0 (L)^{\ast})$ is a covering of degree 6 ramified over a curve $B \subset \Pe^{\ast}_2 $ of degree 18. 
 Our aim is to understand the curves $V$ and $B$ for a general abelian surface $(A,L)$. In this section we study first a special case.
\medskip

\noindent
Let $E$ be an elliptic curve and consider the abelian surface $A = E \times E$. Let $L$ denote the line bundle
$$
L = \co_A (E \times \{ 0 \} + \{ 0 \} \times E + \widetilde{\triangle} )
$$
where $\widetilde{\triangle} = \{ (x , - x) \in E \times E \}$ denote the antidiagonal. According to 
[BL1] $L$ defines an irreducible polarization of type $(1,3)$. Moreover, in [BL2] Proposition 3.3  it is shown that the branch divisor $B \subset \Pe^{\ast}_2$ of $\vp_L : A \ra \Pe^{\ast}_2$ is
$$
B = 3 D \eqno(1)
$$
where $D$ is a plane sextic with with 9 cusps, the dual of the elliptic curve $E$, considered as a plane cubic in $\Pe_2$
embedded by $|3\cdot 0|$ (see [BL2]). We first determine the ramification divisor $R \subset A$ of $\vp_L$.
\medskip

\noindent
We may assume that $L$ is a symmetric line bundle. Then the extended Heisenberg group $H(L)^e$ acts on $B$. We 
may choose the coordinates $(x_0 : x_1 : x_2)$ of $\Pe^{\ast}_2$ in such a way that $H(L)^e$ is generated by (see [CAV])
$$
\sigma = \left(
\begin{array}{ccc}
0&1&0\\
0&0&1\\
1&0&0
\end{array}
\right), \,\, \tau = \left(
\begin{array}{ccc}
1&&\\
& \rho &\\
&& \rho^2
\end{array}
\right) \quad \mbox{and} \quad \iota = \left(
\begin{array}{ccc}
1&0&0\\
0&0&1\\
0&1&0
\end{array}
\right)
$$
where $\rho = \exp (\frac{2 \pi i}{3})$. Consider the  set
$$
{\rm Fix} H(L)^e = \{ 
x \in \Pe^{\ast}_2 \, | \, g (x) = x \ \ {\rm for \ some}\ \  1 \not= g \in H(L)^e \}
$$
of points in $\Pe_2^\ast$ having non trivial stabilizer. It consists of 9 lines and 9 points
$$
\begin{array}{lll}
l_{1i} = \{ x_1 = \rho^i x_0 \}, & \quad l_{2i} = \{ x_2 = \rho^i x_1 \},& \quad l_{3i} = \{ x_0 = \rho^i
x_2 \}\\
y_{1i} = (1 : - \rho^i: 0), &\quad y_{2i} = (0 :1 : - \rho^i), &\quad y_{3i} = (- \rho^i : 0 : 1)
\end{array}
$$
with $i = 0,1,2$. Moreover, the 9 lines intersect in triples in 12 points $x_1, \ldots , x_{12}$, say.
\smallskip

\noindent
This immediately gives that if $x \in l_i, \,\, x \not= x_j$ for all $j$ or if $x = y_i$ for some $i$,
then the orbit of $x$ consists of 9 elements. The orbit of $x_i \,\, (i = 1, \ldots , 12)$ consists of
3 elements. To be more precise, there are exactly 4 orbits of order 3, namely
$$
\begin{array}{lcl}
O(1:0:0)& = & \{ (1:0:0),\quad (0:1:0), \quad(0:0:1) \} \\[1ex]
O(1:1:1)& = & \{ (1:1:1), \quad(1 : \rho : \rho^2),\quad (1: \rho^2 : \rho) \}\\[1ex]
O(1:1:\rho)& = & \{ (1:1:\rho), \quad(1: \rho:1), \quad(\rho:1:1) \}\\[1ex]
O(1:1:\rho^2)& =& \{ (1:1:\rho^2),\quad (1: \rho^2 :1),\quad (\rho^2 :1 :1) \}
\end{array}
$$
All other orbits of $H(L)^e$ in $\Pe^{\ast}_2$ consist of 18 points.
\medskip

\noindent
According to [BL2] Proposition 1.1 with respect to the coordinates $y_0 = 3x^2_0 - 3 \lambda x_1 x_2 ,
\quad y_1 = 3x^2_1 - 3 \lambda x_0 x_2,
\quad y_2 = 3 x^2_2 - 3 \lambda x_0 x_1$ the sextic $D$ is given by the equation
$$
\begin{array}{l}
(y^6_0 + y^6_1 + y^6_2) + 2(2 \lambda^3 -1)(y^3_0 y^3_1 + y^3_0 y^3_2 + y^3_1 y^3_2)\\[1ex]
\quad -6 \lambda^2 y_1 y_2 y_3 (y^3_0 + y^3_1 + y^3_2) - 3 \lambda (\lambda^3 - 4) y^2_0 y^2_1 y^2_2 =0.
\end{array}
$$
Here $\lambda \in \Q - \{1, \rho, \rho^2\}$ is a parameter depending on the elliptic curve $E$. It can
be explicitly determined from the $j$-invariant of $E$. Now an immediate computation using the
above equation gives
\medskip

\noindent
{\bf Lemma 4.1.} {\it For a general elliptic curve $E$ the branch curve $B$ does not contain an orbit of
order $3$ under the extended Heisenberg group} $H(L)^e$.
\medskip

\noindent
Denote by $\triangle := \{ (x,x) \in E \times E \}$ the diagonal, $\Gamma_{-2} = \{ (x, -2x) \in E
\times E \}$ the graph of $(-2)_E$ 
and $\Gamma^t_{-2} = \{ (-2x,x) \in E \times E \}$ its transpose. Then we have
\medskip

\noindent
{\bf Proposition 4.2.} (a): {\it The ramification divisor of $\vp_L: E \times E  \ra \Pe^{\ast}_2$ is 
$$
R = \triangle + \Gamma_{-2} + \Gamma^t_{-2}.
$$
{\rm (b)} $\vp_L | \triangle : \triangle \ra D$ coincides  with the duality map $\triangle = E \ra
E^{\ast}, \,\, x \mapsto {\mathrm {tangent\,\, at}}\,\, x$. In particular $\vp_L | \triangle$ is bijective.}\\
(c) $p_a (R) = 28$.
\medskip

\noindent
{\it Proof.} The map $\vp_L : E \times E \ra \Pe^{\ast}_2$ is a Galois covering with Galois group $D_3$
generated by
$$
T = \left(
\begin{array}{cc}
0 & 1\\
-1& -1
\end{array}
\right) \quad \mbox{and} \quad J = \left( \begin{array}{cc}
0 & 1\\
1 & 0
\end{array}
\right)
$$
(see [BL2] Corollary 3.2). Since $ J | \triangle = \id$, the map $\vp_L$ is ramified in $\triangle$.
Hence
$$
\triangle + \Gamma_{-2} + \Gamma^t_{-2} \subseteq R
$$
since $T(\triangle) = \Gamma_{-2}$ and  $T^2(\triangle) = \Gamma^t_{-2}$.\\
On the other hand $\vp_L | \triangle$ is given by a $2$-dimensional sublinear system of
$$
\begin{array}{lcl}
L | \triangle &=& \co_{\triangle}((E \times \{ 0 \} + \{ 0 \} \times E +
\widetilde{\triangle})\,|\,\triangle)\\[1ex]
&=& \co_{\triangle}(3(0,0) + (\alpha, \alpha) + (\beta, \beta) + (\gamma, \gamma))\\[1ex]
&&\mbox{where} \quad \{0, \alpha, \beta, \gamma\} \quad \mbox{are the division points of} \quad
\triangle.\\[1ex]
&=& \co_{\triangle} (6(0,0)), \quad \mbox{since} \quad \alpha + \beta + \delta \sim 3 \cdot 0 \,\,
\mbox{in} \,\, \triangle.
\end{array}
$$
As $\vp_L (\triangle) =D$ is just the dual curve of $\triangle =E$, this implies (b) and hence (a),
since $B = 3D$. As for (c), $R^2 = 54$ implies $p_a (R) = 28$. \hfill $\square$
\bigskip

\noindent
Using this one can show
\medskip

\noindent
{\bf Proposition 4.3.} {\it Let $V \subset \Pe_2$ denote the discriminant curve of  
$(E \times E, \co_{E
\times E} (E \times \{ 0 \} + \{ 0 \} E + \widetilde{\triangle}))$.
Then
$$
V = 3 E + l_1 + \ldots + l_9
$$
where $E$ denotes the elliptic curve as a plane cubic as above and $l_1, \ldots , l_9$ are the flex tangents of}
$E$.
\medskip

\noindent
{\it Proof.} The projections $\pi_1$ and $\pi_2$ of $A \times \Pe_2$ and the map $\vp_L : A \ra
\Pe_2^{\ast}$ restrict to the following diagram
$$
\begin{array}{ccc}
& \Gamma&\\
 {\scriptstyle \pi_1}\swarrow && \searrow {\scriptstyle \pi_2} \\
R &&V\\
{\scriptstyle \vp} \downarrow \,\,\, &&\\
B&=3D&
\end{array}
$$
where $D$ is the plane sextic with 9 cusps dual to $E$. Let
$$
V = V_0 + l_1 + \ldots + l_v
$$
with lines $l_1 , \ldots , l_v$ and $V_0$ containing no line. 
\smallskip

\noindent
We claim that $V$ contains exactly 9 lines and these are the flex tangents of red $V_0 =E$. To see this
note that for any line $k \subset \Pe^{\ast}_2$ passing through  a cusp of $D$ the branch divisor
of $\vp | \vp^{-1} (k) : C = \vp^{-1} (k) \ra k$ contains a point of multiplicity $\geq 6$, hence $C$ is
singular 
(because $\vp_L$ has degree $6$).
But the pencils of lines through the cusps of $D$ correspond just to the flex tangents $l_1 + \ldots + l_9$ of the
dual curve.\\
On the other hand according to Propositions 3.2 and 4.2 (c) the curves $\Gamma$ and $R$ have the same
arithmetic genus 28. Hence $\pi_1 : \Gamma \ra R$ is an isomorphism apart from possibly contracting
lines of $\Gamma$ to a point in $R$. Applying Proposition 4.2 (a) we conclude that $\Gamma$ contains 3
copies of the elliptic curve $E$ which map under $\pi_2$ to $V$. This implies 
$3 E + l_1 + \ldots + l_9 \subseteq V$. 

It is not difficult to describe the 1-dimensional family of singular curves 
parametrized by $3E\subset V$. The surface $A$ contains a 3-dimensional family 
of singular curves of the form:
$$
E\times\{a\}+\{b\}\times E + t_{(c,c)}^\ast\widetilde{\triangle}
$$
with $a,b,c \in E$. The general curve of this family is easily
seen to  have 3 nodes. Intersecting with $|L|$ we get a 1-dimensional
family of singular curves in $|L|$ the general one having 3 nodes. 
This family is parametrized by the plane cubic $3E\subset V$. 

Consider now a general pencil $\Pe_1 \subset |L|$. The computation made in 
the proof of Prop. 3.4  gives in this case a contribution of 3 from each
curve belonging to $3E$, because each such curve is trinodal. But then since
the number we obtain is 18, and this is also the degree of
$3 E + l_1 + \ldots + l_9 $ we conclude that $V = 3 E + l_1 + \ldots + l_9$
as asserted.
 \hfill $\square$

\medskip\noindent
{\bf Remark 4.4.}  The proof of Prop. 4.3 shows that 
$ \overline{V}_{L,1} =  l_1 + \ldots + l_9$; in particular 
$V \neq \overline{V}_{L,1}$ in this case. In fact the proof shows
that the curves  of $\Pe_1\cap (l_1 + \ldots + l_9)$ contribute 
with multiplicity
one in the computation made in Prop. 3.4, and this can only be possible if these curves are 1-nodal. 

Note that the curve parametrized by a general point of $3E$ is not in  
 $ V_{L,1}$, but it is
 nevertheless a nodal curve (it is in $V_{L,3}$).  In particular the general curve of 
each component of $V$ is a nodal curve.

Finally we remark that from the analysis of this section it follows that the incidence curve
$\Gamma$ of the special surface $(E\times E,L)$ considered is reduced.

%\newpage
\vspace{2cm}
 \begin{center}
{\bf \S 5 The general abelian surface of type  \boldmath$(1,3)$\unboldmath.}
\end{center}
Let $(A,L)$ be a general abelian surface of type $(1,3)$. 
Let $V\subset \Pe_2=|L|$ and 
$B \subset
\Pe^{\ast}_2$ again denote the discriminant curve of $L$ and the branch curve of $\vp_L : A \ra
\Pe^{\ast}_2$. After showing that $V$ and $B$ are irreducible and reduced, we compute their
singularities.

\medskip

\noindent
{\bf Proposition 5.1.} {\it For a general abelian surface $(A,L)$ of type $(1,3)$ we have 
$ V = \overline{V}_{L,1}$ and $V$ is reduced.}
\medskip

\noindent 
{\it Proof.} Since $(A,L)$ is general we may assume that every element of
$|L|$ is irreducible. 
Every curve
belonging to $V$ is irreducible of geometric genus $\leq 3$
because  the genus of the nonsingular elements of $|L|$ is 4. 
Since there are only finitely many curves of geometric genus 2 in $|L|$ (Remark 2.3) 
we deduce that each irreducible component of $V$ consists generically
of curves of geometric genus 3, which therefore either have one node or
one ordinary cusp, and no other singularity. 
Let's consider a family of polarized abelian
surfaces of type $(1,3)$ with general fibre $(A,L)$ and special fibre 
$(E \times E, L_0)$ with $L_0 =
\co_{E \times E} (E \times \{ 0 \} + \{ 0 \} \times E + \widetilde{\triangle})$. 
Under this degeneration  the curve $V$ specializes to a curve $V'$ contained
in the discriminant curve of $|L_0|$. 
A component of $V$ consisting generically of cuspidal curves would degenerate to
a component of $V'$ having the same property.
But we have seen 
(Remark 4.4) that a general element of any component of this curve 
parametrizes nodal curves. So each irreducible component of $V$ consists generically of 
1-nodal curves, i.e. $ V = \overline{V}_{L,1}$.

To prove the last assertion note that under the above degeneration the incidence curve
$\Gamma$ of $(A,L)$ degenerates to the incidence curve of $(E \times E, L_0)$, which is
reduced (Remark 4.4). Therefore $\Gamma$ is reduced as well. The conclusion now follows 
from Proposition 3.4. 
\hfill $\square$

\medskip
\noindent
{\bf Proposition 5.2.} {\it For a general abelian surface $(A,L)$ the branch divisor $B \subset
\Pe^{\ast}_2$ is a reduced and irreducible curve of degree $18$}.
\medskip

\noindent
{\it Proof.} Since $\vp_L (R) =B$ it suffices to show that the  ramification divisor $R \subset A$ is
reduced and irreducible and that $\vp_L :R \ra B$ is birational. As in the proof of
Proposition 5.1 consider a family of polarized abelian
surfaces of type $(1,3)$ with general fibre $(A,L)$ and special fibre $(E \times E, L_0)$ with $L_0 =
\co_{E \times E} (E \times \{ 0 \} + \{ 0 \} \times E + \widetilde{\triangle})$. Certainly here the
ramification divisor $R$ of $\vp : A \ra \Pe^{\ast}_2$ specializes to the ramification divisor $R_0$ 
of $\vp_{L_0} : E \times E \ra \Pe^{\ast}_2$. According to Proposition 4.2 (a) 
$$
R_0 = \triangle + \Gamma_{-2} + \Gamma^t_{-2} .
$$
Now $R^2 = R^2_0 = 54$. This implies that $\co_A (R)$ is of type $(1,27)$ or $(3,9)$. The abelian
surface $(A,L)$ being general, it cannot be of type $(1,27)$. Hence
$$
\co_A (R) \equiv L^3.
$$
If $R$ were  not irreducible and reduced, say $R = R_1 + R_2$ or $R_1 + R_2 + R_3$, then we would
have $\co_A (R_1) \equiv L$. But $R_1$ specializes to one of the curves $\triangle, \Gamma_{-2}$ or
$\Gamma^t_{-2}$, say to $\triangle$. This would give $\co_{E \times E} (\triangle) \equiv L_0$, a
contradiction, since $L^2_0 =6$ and $\triangle^2 =0$.\\
It remains to be proved that $\vp_L : R \ra B$ is  birational. 
For any point $p \in \Pe^{\ast}_2$ we have that $\vp^{-1}(p)$ consists of the $6$ base points
(counted with multiplicity) of a pencil in $|L|$ which is the  line of  $\Pe_2$ dual to $p$. 
Since $B=\vp(R)=\vp(\pi_1(\Gamma))$, we have that $p\in B$ if and only if the corresponding
pencil is of the form $|L(-x)|$ for some $(x,[C]) \in \Gamma$, (and of course $\vp(x)=p$). 
If $p$ is a general point of
a component of $B$ then  $[C]\in V_{L,1}$ and $x$ is the node of $C$ (by Prop. (5.1));
moreover  the pencil $|L(-x)|$ is the tangent
line to $V_{L,1}$ at $[C]$ (see \S 1).
Assume that $\vp_L : R \ra B$ is not birational. 
Then for a general choice of $p\in B$ there are $(x_1,[C_1]), (x_2,[C_2])  \in \Gamma$ 
such that $p= \vp(x_1) = \vp(x_2)$ and $x_1 \ne x_2$.  
If $C_1=C_2=C$ then the curve $C$ has two nodes 
and therefore it is not in $V_{L,1}$, a contradiction. If $C_1 \ne C_2$ then the line 
$|L(-x_1)|=|L(-x_2)|$ of $\Pe_2$ is bitangent to $V=\overline{V}_{L,1}$ 
at the points $[C_1]$ and $[C_2]$. 
This contradicts the generality of $p$ because the reduced curve $V$ 
has only finitely many bitangent lines.
 This completes the proof of Proposition 5.2. \hfill $\square$
\medskip

\noindent

{\bf Proposition 5.3.} 
{\it The  curve $V$ of a general abelian surface $(A,L)$ of type $(1,3)$ is reduced of degree 18 and 
admits at most nodes, cusps or ordinary tacnodes as singularities.}
\medskip

\noindent
{\it Proof.} 
By 5.1 and 3.4   $V$ is reduced of degree 18.
If $[X] \in V$ then, thanks to Proposition 2.2, one of the following occurs:

i) $X$ has one node ($g=3$)

ii) $X$ has one cusp ($g=3$)

iii) $X$ has two nodes ($g=2$)

\noindent
and there are no other possibilities. 

Case i) occurs if and only if $[X] \in V_{L,1}$. In this case $[X]$ is a nonsingular point 
of $V$. 

Suppose we are in case ii). The natural map 
$h: \frac{H^0(L)}{\langle\sigma_X\rangle} \ra H^0(T^1_X)$ (see \S 1)
is bijective: in fact domain and codomain are both of dimension 2 and $Im(h)$ contains two 
independent vectors corresponding to an infinitesimal deformations which smooths $X$
and to an infinitesimal deformation which deforms $X$ to a nodal curve. Therefore $V$ induces,
locally around $[X]$, a semiuniversal deformation of the cusp. It is well known that the 
locus parametrizing singular fibres is an ordinary cusp at $[X]$ (see e.g.
[D-H], p. 3).

Suppose finally that we are in case iii). 
By taking a general pencil $\Pe_1 \subset |L|$ containing $[X]$ and applying 
the same argument of Proposition 3.4 we see that $[X]$ absorbs two of the 18 intersections
$\Pe_1 \cap V$, because $e(X)=e(X_{gen})+2$. Therefore $[X]$ is a double point of $V$.
Let $N_1, N_2$ be the nodes of $X$. For each
$i=1,2$  the vector space 
$$
\frac {H^0({\cal I}_{N_i}\otimes L)}{\langle \sigma_X \rangle}
\subset \frac {H^0(L)} {\langle \sigma_X \rangle}
$$
is the tangent space to the local analytic family of deformations of $X$ which keep the
node $N_i$. Each such local family is a linear branch of $V$ at $[X]$ (see [D-H] again). 
Therefore $V$ has two
linear branches at $[X]$, which is therefore either an ordinary  node or an ordinary tacnode.
The last possibility occurs precisely when 
$H^0({\cal I}_{N_1}\otimes L) = H^0({\cal I}_{N_2}\otimes L)$, and this happens 
if and only if $N_1$ and $N_2$ are in the same fibre of $\vp_L$. 
\hfill $\square$

\medskip\noindent
{\bf Proposition 5.4.}   {\it Let $(A,L)$ be a general abelian surface  of type $(1,3)$ and let 
$V = V'+l_1+\cdots +l_v$, where $l_1,\ldots,l_v$ are the 
line components of $V$. Then $B$ and $V'$ are dual to each other
(in particular $V'$ is irreducible).}

\medskip\noindent
{\it Proof.}  As seen in the proof of Prop. 5.2, 
for any point $p \in \Pe^{\ast}_2$ we have that $\vp ^{-1}(p)$ consists of the $6$ base points
(counted with multiplicity) of a pencil in $|L|$ which is the  line of  $\Pe_2$ dual to $p$ 
if $p$ is a general point of
a component of $B$ the curve $C$ is nodal at $x$, $[C]\in V'$ and the pencil $|L(-x)|$ is the tangent
line to $V'$ at $[C]$. This proves that $B=V'^\vee$,  the dual curve of $V'$. The converse follows from
the fact that $V'^{\vee\vee}=V'$.\hfill $\square$

\medskip

\noindent
Recall that the projections $\pi_1$, and $\pi_2$ and the map $\vp_L$ restrict to the following diagram
$$
\begin{array}{ccc}
& \Gamma &\\
 {\scriptstyle \pi_1} \swarrow && \searrow {\scriptstyle \pi_2} \\
R && V\\
{\scriptstyle \vp_L} \downarrow \,\,\, &&\\
B &&
\end{array}
$$
Let $l_1, \ldots , l_v$ denote the lines in $V$ and $V = V' + l_1 + \ldots + l_v$. According to
Proposition 5.4 the curve $V'$ is the dual of $B$, and it  is irreducible. The
lines $l_i$ correspond to points on $R$. According to Lemma 4.1 we may assume that the number $v$ of
lines is divisible by $9$, since $(A,L)$ is general. As $B$ is dual to $V'$, we are left with 2 cases .
\begin{itemize}
\item[I.] $V = V'$ irreducible  of degree $18$ in $\Pe_2$ and $v =0$
\item[II.] $V = V' + l_1 + \ldots + l_9$ with $V'$ irreducible of degree 9 in $\Pe_2$.
\end{itemize}
For the curve $\Gamma$ this implies: Either $\Gamma = \Gamma'$ is irreducible or $\Gamma = \Gamma' +
l'_1 + \ldots + l'_9$ with $\Gamma'$ irreducible and $\pi_1 (l'_i)$ is a point in $R$ for $i = 1, \ldots
, 9$.
\medskip

\noindent
{\bf Proposition 5.5.} {\it The curve $\Gamma'$ is smooth.}
\smallskip

\noindent
{\it Proof.} By Proposition 5.3   $V'$ and hence also  $\Gamma'$ admit at most nodes, cusps and ordinary
tacnodes.
\smallskip

\noindent
If $\Gamma'$ has a node or an ordinary tacnode in $(x, [C])$, Then $V'$ also has a node or a tacnode in
$[C]$. This implies that the curve $C$ is of genus 2. But $C$ admits at most nodes as singularities by
Proposition 2.2. Hence $\pi^{-1}_2 ([C])$ consists of 2 points, so both of them have to be smooth, a
contradiction.
\smallskip

\noindent
If $\Gamma'$ has a cusp in $(x,[L])$, then $V'$ has a cusp in $[C]$. But $\Gamma'$ projects birationally onto
the dual curve $B$ of $V'$ and the dual of a cusp is a smooth point. Hence $\vp_L (x)$ is a smooth point
of $B$ and so is $(x,[C])$ in $\Gamma'$, since all maps are birational. \hfill $\square$
\medskip

\noindent
Now we are in a position to prove the main theorem of this paper.
\medskip

\noindent
{\bf Theorem 5.6.} {\it Let $(A,L)$ be a general abelian surface of type $(1,3), \,\,\,
B \subset \Pe^{\ast}_2$ the branch locus of $\vp_L : A \ra \Pe^{\ast}_2$ and $V \subset \Pe_2 = |L|$ the 
closure of the Severi variety  $V_{L,1}$.
\smallskip

\noindent
Both $V$ and $B$ are irreducible of degree $18$ in $\Pe_2$ (resp. $\Pe^{\ast}_2$) smooth apart from $72$ 
cusps and $36$ nodes or tacnodes (each tacnode counting as two nodes), and admit $72$ flexes and $36$ bitangents.}
\medskip

\noindent
{\it Proof.}  By Proposition 5.3 \,\,\,
$V$ admits at most nodes, cusps or ordinary tacnodes as singularities. Since for the following computations an ordinary 
tacnode counts exactly the same as $2$ nodes, we assume for simplicity that $V$ admits at most nodes and cusps.
\smallskip

\noindent
Let $V'$ the irreducible component of $V$ as above. Denote by $\nu$ (respectively 
 $\kappa, f,$ and $b$) the number of nodes (respectively cusps, flexes and bitangents) 
 of $V'$. The classical Pl\"{u}cker  formulas (see [Wa]) say that these numbers are related to the degree $d$ 
of $V'$ and the degree $m$ of its dual plane curve $B$ as follows:
$$
\begin{array}{lr}
m = d(d-1) - 2 \nu - 3 \kappa & \hspace*{6cm}{\rm (1)}\\
d = m(m-1) -2b -3f & \hspace*{6cm}{\rm (2)}\\
f = 3d(d-2) - 6 \nu - 8  \kappa & \hspace*{6cm}{\rm (3)}
\end{array}
$$
Suppose first that we are in case II: \,\, $V = V' + l_1 + \ldots + l_9$. Let $\Gamma = \Gamma' + l'_1 +
\ldots l'_9$ denote the corresponding decomposition. Then $\Gamma'$ is smooth by Proposition 5.5 and by
Proposition 3.2   we have
$$
p_a (\Gamma) = p_a (R) = 28.
$$
For this situation there are only 2 possible cases:
\begin{itemize}
\item[(i)] The lines $l_1 , \ldots , l_9$ intersect $\Gamma'$ transversally in 1 point each and $\pi_2 :
\Gamma' \ra R$ is an isomorphism.
\item[(ii)] The line $l_i$ intersects $\Gamma'$ transversally in 2 points $x_i$ and $x'_i$ for $ i = 1,
\ldots, 9$ and $\pi_2 (x_i) = \pi_2 (x'_i)$ is a node in $R$ whereas $\pi_2 : \Gamma' \ra R$ is an
isomorphism elsewhere (here we use the fact that the curve $R$ is Heisenberg invariant and that each
$H^0(L)^e$-orbit consists of $\geq 9$ points).
\end{itemize}
Suppose first we are in case (i): Then $g(V') = g(\Gamma') = g(R) = 28$. On the other hand $V'$ is a
plane curve of degree $9$. Hence
$$
28 = p_a (V') = g (V') + \nu +\kappa  = 28 + \nu + \kappa,
$$
so $\nu  = \kappa = 0$. But then (1) says $18 = 72$, a contradiction. In the second case $g(V') = g(\Gamma) =
p_a(R) -9 =19$. Hence 
$$
28 = p_a (V') = 19 +\nu + \kappa .
$$
So $\nu +\kappa = 9$. On the other hand (1) says
$$
2 \nu = 3 \kappa = 54.
$$
This implies $\kappa  = 18, \nu =  27$, a contradiction.
\smallskip

\noindent
Hence $V = V'$ is irreducible of degree 18 in $\Pe_2$. So $\Gamma$ is the normalization of $V$ and $g(Y)
= g(\Gamma) = 28$ and we have
$$
136 = p_a (V) = g (V) + \nu + \kappa = 28 + \nu +\kappa
$$
i. e.
$$
\nu + \kappa  = 108.
$$
On the other hand (1) says
$$
2 \nu + 3 \kappa  = 288.
$$
Combining both equations we obtain $\kappa  = 72$ and $\nu = 36$. But then (2) and (3) give $f = 72 $ and $b 
= 36$. 

By Proposition   5.4   the statement follows also for  $B$. \hfill $\square$

\medskip 
\noindent
{\bf Remark.} The 36 nodes of $V$ correspond to the curves of genus 2 in $|L|$, which are all 2-nodal.
The 72 cusps correspond to the curves of genus 3 in $|L|$ having a cusp.
We don't know whether tacnodes actually occur for general $(A,L)$
(actually for any $(A,L)$).

\newpage
\noindent
\begin{center}
{\bf References}
\end{center}
\begin{itemize}
\item[{[BL1]}] Ch. Birkenhake, H. Lange: Moduli Spaces of Abelian Surfaces with Isogeny.
Proc. of the Conference ``Geometry and Analysis'', TATA Institute Bombay, Oxford
University Press (1995), 225-243.
\item[{[BL2]}] Ch. Birkenhake, H. Lange: A family of abelian surfaces and curves of genus 4, manuscr.
math. 85 (1994), 393 - 407.
\item[{[BPV]}] W. Barth, C. Peters, A. van de Ven: Compact complex surfaces. Erg. d. Math. Springer,
Berlin (1984).
\item[{[CAV]}] H. Lange, Ch. Birkenhake: Complex Abelian  Varieties. Grundlehren Math. Wiss. 302.
Berlin-New York 1992.
\item[{[D-H]}] S. Diaz, J. Harris:  Geometry of the Severi variety I, Trans. AMS 309 (1988), 1 - 34.
\item[{[W]}] J. Wahl: Deformations of plane curves with nodes and cusps, 
Am. J.  Math. 96 (1974), 529-577.
\item[{[Wa]}] R. J. Walker: Algebraic Curves, Princeton University Press, 1950.
\end{itemize}

\vskip1cm\noindent
Mathematisches Institut
\hfill\break
Bismarckstr. $1{1\over 2}$, D-91054 Erlangen, Germany

\bigskip\noindent
Dipartimento di Matematica, Universit\`a Roma III
\hfill\break
L.go S.L. Murialdo 1, 00146 Roma (Italy)
\end{document}